\documentclass{ambphack3}

\usepackage{amsfonts}

\begin{document}

\newcommand{\lie}{{\rm Lie}}
\newcommand{\ssi}{S^*S^1}
\newcommand{\cons}{\frac{\Gamma\left(\frac{n}{q}\right){\rm
dim}(E)}{q^2(2\pi)^n }}
\newcommand{\ext}{{\rm ext}}
\newcommand{\id}{{\rm id}}
\newcommand{\ev}{{\rm ev}}
\newcommand{\Iso}{{\rm Iso}}
\newcommand{\ord}{{\rm ord}}
\newcommand{\MMaps}{{\rm Maps}}
\newcommand{\Fr}{{\rm Fr}}
\newcommand{\BbP}{{\mathbb P}}
\newcommand{\BbC}{{\mathbb C}}
\newcommand{\BbZ}{{\mathbb Z}}
\newcommand{\pdo}{\Psi{\rm DO}}
\newcommand{\Tr}{{\rm Tr}}
\newcommand{\tr}{{\rm tr}}
\newcommand{\sm}{S^*M}
\newcommand{\res[1]}{{\rm res_w}#1}
\newcommand{\ress}{{\rm res_w}}
\newcommand{\cl[1]}{{\mathcal Cl}_{#1}}
\newcommand{\calAA}{\mathcal A}
\newcommand{\calA}{{\rm Ell}({\mathcal E})}
\newcommand{\calB}{{\mathcal Cl}({\mathcal E})}
\newcommand{\calD}{{\mathcal D}}
\newcommand{\calE}{{\mathcal E}}
\newcommand{\calG}{{\mathcal G}}
\newcommand{\adP}{{\rm ad}\ P}
\newcommand{\AdP}{{\rm Ad}\ P}
\newcommand{\R}{\mathbb R}
\newcommand{\Z}{\mathbb Z}
\newcommand{\C}{\mathbb C}
\newcommand{\N}{\mathbb N}
\newcommand{\e}{\varepsilon}
\newcommand{\Aut}{{\rm Aut}}
\newcommand{\End}{{\rm End}}
\newcommand{\Hom}{{\rm Hom}}
\newcommand{\ad}{{\rm ad}}
\newcommand{\Ad}{{\rm Ad}}
\newcommand{\calC}{{\mathcal C}}
\newcommand{\dvol}{{\rm dvol}}
\newcommand{\oo}{{\rm o}}
\newcommand{\OO}{{\rm O}}
\newcommand{\asy}{{\rm asy}}
\newcommand{\fp}{{\rm f.p.}_{\e=0}}
\newcommand{\Ell}{{\mathcal Ell}}
\newcommand{\PP}{\mathbb P}

\title{Chern-Weil Constructions on $\pdo$ Bundles}

\FirstAuthor{Sylvie Paycha}
\FirstDepartment{Laboratoire de Math\'ematiques}
 \LongFirstDepartment{Appliqu\'ees}
\FirstInstitution{Universit\'e Blaise Pascal}
\LongFirstInstitution{(Clermont II)}
\FirstStreetAddress{Complexe Universitaire}
\LongFirstStreet{de C\'ezeaux}
\FirstCity{63177 Aubi\'ere Cedex}
\FirstCountry{France}
\FirstEmail{paycha@ucfma.univ-bpclermont.fr}
\SecondAuthor{Steven Rosenberg}
\SecondDepartment{Department of Mathematics}
\LongSecondDepartment{and Statistics}
\SecondInstitution{Boston University}
\SecondStreetAddress{111 Cummington St.} 
\SecondCity{Boston, MA 02215}
\SecondCountry{USA}
\SecondSupport{Partially supported by the NSF and the CNRS.}
\SecondEmail{sr@math.bu.edu}

\maketitle 

\ShortAuthor{S. Paycha and S. Rosenberg}

\ShortTitle{Characteristic Classes on $\pdo$ bundles}


\begin{abstract}  We construct  Chern-Weil classes 
on infinite dimensional vector bundles with structure group contained
in the algebra $\cl[\leq 0](M, E)$ of non-positive order classical
pseudo-differential operators acting on  a finite rank
vector bundle $E$ over a closed manifold $M$.  Mimicking
the finite dimensional Chern-Weil construction, we replace the ordinary
trace on matrices by linear functionals on $\cl[\leq 0]
(M, E)$ built from the leading symbols of the operators.
The corresponding Chern classes vanish on loop groups, but a weighted
trace construction yields a non-zero class previously constructed by Freed.
For loop spaces, 
the structure group reduces to a gauge group of
bundle automorphisms, and
we produce
non-vanishing universal Chern classes in all degrees, using a universal
connection theorem for these bundles. 
\end{abstract}


\section{Introduction}

Infinite dimensional vector bundles with connections
are frequently encountered in
mathematical physics; basic examples include the tangent bundle to
loop spaces and the bundle of spinor fields associated to a family of
Dirac operators.  In this paper, we construct  Chern forms
and Chern classes for a class of vector bundles including the tangent
bundle to loop spaces, and produce examples of non-vanishing
classes.
Motivation for this work comes
from \cite{B, BF,CDMP,  jp, PR}, among others.

We should emphasize that we do not construct a unified Chern-Weil
theory as free from choices as the finite dimensional theory.  Some
choices seem inevitable: different choices of model spaces for the
vector bundle fibers lead to inequivalent topologies, and the
natural choice of the general linear group of a Hilbert space as
structure group leads to a trivial
theory by \cite{Ku}.  
Other choices seem more arbitrary: for the classical groups
$GL(n,\C), U(n)$ and (most of) $SO(n)$, Chern-Weil theory depends
crucially on properties of the ordinary trace on matrices.  In
infinite dimensions,
there are many
possible generalizations of the trace, the operator trace being the
most obvious. (It is seldom applicable, since the curvature
of the connections we consider are not trace class operators.)
Much of the paper is devoted to isolating the key features of the
trace, and then constructing useful generalizations to infinite
dimensions.  

Motivated by the examples mentioned, we concentrate on
bundles $\calE\to B$ 
with fibers $\calE_b$ modeled on either Sobolev spaces of sections or
smooth sections of a finite rank bundle $E$ over a closed manifold
$M$.  The structure 
group is a group of classical pseudo-differential operators
($\pdo$s) on $E$, and the curvatures of the connections we consider
take values in a corresponding Lie algebra of $\pdo$s.  For loop
spaces, the structure group is typically multiplication operators, but
the curvatures of the connections we study take values in $\pdo$s.

In more detail, in \S2 we discuss traces as morphisms $\lambda:\AdP\to
\C$ from the adjoint bundle of a principal bundle $P$ to the trivial $\C$
bundle.  Here the structure group and the base may be infinite
dimensional, and we are thinking of $P$ as the principal bundle
associated to $\calE.$
  We show that such a functional acting on the curvature
$\Omega$ of a connection on $P$
leads to characteristic
forms and classes on the base generalizing the usual forms
$\tr(\Omega^k)$ in finite dimensions (Theorem 2.2).  

In \S3, we introduce two types of traces in infinite
dimensions. The first type is the Wodzicki residue, the unique trace
on the space of classical $\pdo$s.  
However, we can often restrict the structure group to invertible
zeroth order $\pdo$s.  The corresponding Lie algebra of non-positive
order $\pdo$s admits a family of traces of the form $A\mapsto \Lambda(
\sigma_0^A)$, where $\sigma_0^A$ is the leading symbol of $A$ and
$\Lambda$ is any distribution on the cosphere bundle of $M$.  
For both traces, the corresponding Chern classes vanish for loop groups.
We show
that the associated universal
Chern classes are non-zero  for the
structure group of loop spaces, a gauge group
 (Theorem 3.3).  This result depends on
a universal connection theorem proved in \S4. 
Based on Freed's work \cite{F}, we construct a nontrivial
characteristic class on loop groups in \S3.3.  To do so, we
relate the leading term in the asymptotic expansion $\tr(Ae^{-\e
Q})$, where $Q$ is a generalized Laplacian, to a trace
$\Lambda^Q(\sigma_0^A)$, where $\Lambda^Q$ is a distribution depending
on $Q$.  We then
 analyze Freed's  conditional first
Chern form on loop groups as the most divergent term in the asymptotic
expansion of a weighted trace, and give a short, direct proof that
this form is closed

In \S4, we prove that the bundle $E\calG\to B\calG$ has a universal
connection, where $\calG$ is the gauge group of the finite rank bundle
$E$.  We hope this result is of independent interest.

The second author would like to thank the Universit\'e Blaise Pascal
for its hospitality during the preparation of this article.
Conversations with Simon Scott on this subject are also gratefully
acknowledged.

\section{Abstract Chern-Weil calculus}

We first give a quick review of connections on principal
bundles and associated bundles; see \cite{BGV} for more details.
Let $P\to B$ be a Fr\'echet, resp.  
Hilbert  principal $G$-bundle over a manifold $B$, 
with $G$ a Fr\'echet, resp.  Hilbert
 Lie group, let $\rho: G\to \Aut(W)$ be a representation of $G$ on
 some vector space $W$ and let ${\cal W}:= P\times_\rho W\to B$ be the
 associated vector bundle.  
It is a Fr\'echet, resp.  Hilbert  vector bundle over $B$.  

For example,
if $E\to B$ is a finite
 dimensional vector bundle with structure group $GL(n,\C)$, then
the frame bundle $P := GL(E)\to B$ is a principal $GL(n,\C)$-bundle
 and $\AdP$ is isomorphic to $\End(E)$, the endomorphism bundle of
 $E$.  If $E$ is an infinite dimensional vector bundle with typical
fiber $V$ and structure
group $G$, then the associated principal $G$-bundle $P^E$ is built by
gluing copies of $G$ over a point $b\in U_i\cap U_j$ in overlapping
charts by the same $g = g_b\in G$ that glues the copies of $V$ over $b$.

 The space $\Omega(B, {\cal W})$ of ${\cal W}$-valued forms on $B$ can
 be identified with the space \\
$\left( \Omega(P) \otimes
 W\right)_{\rm basic}$
of basic forms   $\alpha$ on $P$ with values in the trivial
 bundle $P\times W$.  Recall that a form is basic if it is $G$-invariant
($\rho(g) (g^* \alpha)= \alpha$) and horizontal
($i_{\bar X} \alpha= \alpha$ for the
canonical vector field $\bar X$ generated by $X\in A = {\rm Lie}(G)$).

   A connection on $P$ given by a one-form $\theta\in \Omega^1(P, A)$
with $\Ad_g (g^* \theta)=\theta$ for $g\in G$ and $i_{\bar X}
\theta=X$ yields a covariant derivative on ${\cal W}$ in the
following way. The derivative $D\rho:A\to \End(W)$ of
$\rho$ induces $D\rho(\theta)\in
\Omega^1(P)\otimes \End(W)$ and hence a connection $\nabla= d+
D\rho(\theta)$ on the trivial bundle $P\times W$.  The resulting
operator $\nabla$ on $\Omega(P)\otimes W $ descends to an operator
on $\left( \Omega(P) \otimes W\right)_{\rm basic} $ and hence a
covariant derivative on ${\cal W}$ also denoted by $\nabla$, with
$\nabla= d+D\rho(\theta)$ locally.  Starting with a vector bundle
$E\to B$ with connection $\nabla$, we produce a connection one-form
$\theta$ on $P^E$, and the curvature of $\nabla$ takes values in
$\AdP.$  

If $W=A$ is the Fr\'echet Lie algebra of
 $G$ and if $\rho$ is the adjoint action, then ${\cal W}$ is called
 the adjoint bundle and denoted by $\AdP$. 
 Recall that the adjoint representation
$\Ad:G\to \Aut(A)$  is the differential of conjugation in $G$: 
$\Ad_ga := (D_e
 C_g) a,$ where $C_g:G\to G$ is $C_g(h) = ghg^{-1}.$
The differential of $\Ad$, $\ad = D\Ad:A\to \End(A)$, is given by 
$\ad_b(a) = [b,a]$, since
$$\ad_b(a)=\frac{d}{dt}\left|_{_{_{_{_{_{_{_{t=0}}}}}}}} \Ad_{e^
{tb}}(a)=\frac{d}{dt}\left|_{_{_{_{_{_{_{_{t=0}}}}}}}}
\frac{d}{ds}\left|_{_{_{_{_{_{_{_{s=0}}}}}}}}
C_{e^{tb}}e^{sa}=[b,a].\right.\right.\right.$$ 
In particular, a connection one-form $\theta\in \Omega^1(P, A)$ yields a
connection $\nabla^\ad $ on $\AdP$, with $\nabla^\ad= d+
[\theta,\cdot]$. (For this reason, our $\AdP$ is often denoted
$\adP.$) If $P = GL(E)$, a connection $\nabla^E$ on $E$
induces a connection one-form on $P$, and the induced connection on
$\AdP$ coincides with the connection $\nabla^{\Hom(E)}$ associated to
$\nabla^E.$  
\medskip

A linear form on $A$:
$$\lambda: A \to  \C$$
such that $\Ad^*\lambda:= \lambda\circ \Ad=\lambda$ induces a bundle
 morphism 
$$\lambda: \AdP \to B\times \C$$ defined as follows. Given  a 
local
trivialization $(U, \Phi)$, where $U\subset B$ is open
and  $\Phi: \AdP|_{U}\to U\times A$ is an isomorphism, and a
 local section $\sigma\in \Gamma(\AdP|_{U})$, we set
$$\lambda (\sigma) := \lambda (\Phi
(\sigma)). $$
This definition is independent of the
local trivialization. Indeed,  given another local
trivialization $(V,\Psi)$, at $b\in U\cap V$ we have
$$\lambda (\Phi
(\sigma))= \lambda (\Ad_g\Psi
(\sigma)) =  \lambda ( \Psi
(\sigma)),\quad \hbox{for some }  \  g = g_b\in G.$$
 
The connection $\nabla^\ad$ on $\AdP$ induced by a connection $\theta$ on
$P$ induces in turn a connection $\nabla^*$ on the dual bundle
$\AdP^*$ (i.e.~$(\AdP)^*$), 
which is locally described by ${\nabla^\ad}^*=d+\ad_{\theta}^*$.
Since $\Ad^*\lambda=\lambda$ implies $\ad^*\lambda=0$, 
we have ${\nabla^\ad}^* \lambda= d\lambda=0$, since
$\lambda$ is locally constant.  Summarizing, we have:

\begin{lemma}\label{mainlemma}
Let $ \lambda: \AdP \to B\times \C$ be the linear morphism induced by a
linear form $\lambda: A\to \C$ with  $\Ad^*
\lambda=\lambda$.  Let $\nabla^\ad= d+ [\theta, \cdot] = d + \ad_\theta$
 be a connection
on $\AdP$ induced by a connection $\theta$ on $P$.  Then
\begin{equation}\label{assumption}
d\circ\lambda=\lambda\circ\nabla^\ad.\end{equation} 
\end{lemma}

\begin{proof} Since $d\lambda=0$, 
we have
$d\circ\lambda= 
\lambda\circ d$ locally. However, $ \ad^* \lambda=0$ implies
$\lambda\circ d =  \lambda\circ (d+\ad_\theta)=\lambda\circ\nabla^\ad,$ so 
$d\circ\lambda=\lambda\circ\nabla^\ad$ globally.
\end{proof}

Abusing notation, we will sometimes denote $\nabla^\ad \alpha $ by
 $[\nabla, \alpha]$, for $\alpha$ an $\AdP$-valued form,
in analogy to the local description $\nabla^\ad = d+ [\theta, \cdot]$, with the
 understanding that $[\nabla,\alpha]$ is a superbracket with respect
 to the $\Z_2$-grading on differential forms.
\medskip 

The lemma gives the main result of this section.

\begin{proposition}\label{maintheorem}  Let $ P = P^E$ be the principal bundle associated to a
vector bundle $E\to B$ with structure group $G$ and connection
$\nabla$ and curvature $\Omega.$
  Assume that the associated connection $\nabla^\ad$ on $\AdP$
has $d \circ \lambda= \lambda\circ \nabla^\ad.$
Then for any analytic
function $f:\C\to\C$, the form $\lambda(f(\Omega) )$ is closed, and its 
de Rham cohomology
class in $H^*(B;\C)$ is independent of the choice of
$\nabla$.
  \end{proposition}

As usual, we mean that the degree $k$ piece of $\lambda(f(\Omega) )$
is a closed $2k$-form, for all $k\in\N.$
\medskip

\begin{proof} The
usual finite dimensional proof (see e.g.~\cite{BGV}) runs through, with
  ordinary traces replaced
by  $\lambda$.

In more detail, $\lambda(f(\Omega))$ is closed because $\lambda (
\Omega^k)$ is closed for any $k\in \N$, which we check in a local
trivialization of $\AdP$.
We have:
$$ d\,\lambda(\Omega^k)
 =
\lambda(\nabla^\ad \Omega^k  )
 = \sum_{j=1}^k \lambda(\Omega^{j-1}( \nabla^\ad
\Omega )  \Omega^{k-j})
= 0$$
where we have used  the Bianchi identity $\nabla\Omega =0$
in the last identity.

To check that the corresponding de Rham class is independent of the choice
of connection, we
consider a smooth one-parameter family of   connections
$\{\nabla_t, t\in \R\}$ on $E$ and the induced 
family of connections $\nabla_t^{\rm ad}$ on $\AdP.$
Then
\begin{eqnarray}\label{cw}
 \frac{d}{dt}\lambda(\Omega_t^k)
&=&
 \sum_{j=1}^k \lambda\left(\Omega_t^{k-j} \left(
\dot \nabla_t \nabla_t + \nabla_t  \dot
\nabla_t\right)
\Omega_t^{j-1}\right)  
=
 \sum_{j=1}^k \lambda\left(\Omega_t^{k-j}
(\nabla_t^\ad
 \dot\nabla_t) \Omega_t^{j-1}\right)
\\
 &=&
 \sum_{j=1}^k \lambda\left(
\nabla_t^\ad\Omega_t^{k-j}  \dot \nabla_t
\Omega_t^{j-1}\right)
=  d  \sum_{j=1}^k
\lambda\left(
\Omega_t^{k-j}
\dot \nabla_t
\Omega_t^{j-1} \right). \nonumber\end{eqnarray}
In the first equality, we use $\nabla_t^2 = \Omega_t$, in the second 
 we have extended the bracket connection to forms, and in the third
 we have used the Bianchi identity.
(\ref{cw}) shows that the dependence on the connection is measured by an
exact form and hence
vanishes in cohomology.\\
${}$  \end{proof}

\begin{remarks} (1) This proof does not use the linearity of
$\lambda$, and extends to any map $\lambda:\AdP\to\C$ satisfying
(2.1).

(2) Just as in finite dimensions, if $E$ is a trivial bundle, then
$[\lambda(f(\Omega))] =0$, because we can connect $\nabla$ to the
trivial connection.  In particular, if $M$  is a
finite dimensional manifold, and $F\to N$ is a trivial finite
rank bundle over a finite dimensional manifold $N$,
there is a
corresponding infinite dimensional  trivial  bundle $E$ over
$B := {\rm Maps}(M,N)$,
whose fiber over $\phi$ is the space
of sections $\Gamma(\phi^*F)$   for $\phi^*F\to M.$  For example, we
can take $F = TN$ (if $N$ is parallelizable), in which case $E = TB.$
Thus characteristic classes formed by such $\lambda$'s will vanish on
$TB$, confirming results in \cite{jp}.  In particular, for $M=S^1$ and
$N=G$, a Lie group, we get the vanishing of the characteristic classes
for the tangent bundle to the loop group $LG.$
\end{remarks}

The previous theorem yields the usual Chern-Weil classes:

\begin{corollary}
Let $ G\subset GL(n,\C) $ be a finite dimensional Lie group, and let
$E\to B$ be a vector bundle with structure group $G$. Let
$\nabla$ be a connection on $E$ with curvature $\Omega.$
For
any analytic function $f$, the forms $\tr(f(\Omega))\in \Omega^*(B,\C)
$ are closed and their de Rham cohomology
 classes are independent of the choice of
$\nabla$. \end{corollary}

This follows from Proposition 2.2 by passing from $E$ to $P^E$
and using
$\lambda= \tr$, the ordinary trace on matrices.

\begin{remark} For $GL(n,\C)$ and $U(n)$, all characteristic classes
are generated by $\tr(\Omega^k), k\in \N.$ However, we do not capture
the Euler class for $SO(n,\R)$ by this procedure, as this class is
generated by the non-linear, but Ad-invariant function $\sqrt{\det}.$
We can treat this case using 
the identity $\det(1 + A) = \sum_k \tr(\Lambda^k A)$, or by the
previous remark (1).  
\end{remark}

\begin{notation} Throughout the paper, ``$\pdo$s'' means
classical pseudo-differential operators, \hfill\\
and 
$\cl[](M,E)$ denotes the
space of all classical $\pdo$s acting on smooth sections of the finite
dimensional Hermitian bundle
$E$ over a closed Riemannian 
manifold $M$. $\cl[k](M,E)$ denotes the subspace of
$\pdo$s of order $k\in\R$, $\cl[\leq k](M,E)$ denotes the space of
$\pdo$s of order at most $k$. $\cl[k]^*(M,E)$ denotes the set of
invertible operators in $\cl[k](M,E)$.  $\Ell_k = 
\Ell_k(M,E)$ denotes the space of
elliptic operators of order $k$, and $\Ell_k^*, \Ell_k^+$ denote the
spaces of invertible, resp.~positive elliptic operators of order $k$. 
\end{notation}

A bundle $\calE$ with fiber modeled on $C^\infty(M,E)$
or on 
$H^s(M,E)$ is a {\it $\Psi$DO bundle} if the transition
maps lie in the Lie group $ \cl[0]^*(M,E)$.  Here $C^\infty(M,E),
H^s(M,E)$ are the spaces of smooth and $s$-Sobolev class
sections of $E$, respectively.
Strictly speaking, the
 transition maps must be zero order only
in the Hilbert (i.e.~Sobolev) setting, since there the transition maps and
their inverses must be bounded invertible operators on a Hilbert
space. In the smooth setting, one can allow non-zero order $\pdo$s.

A connection on a $\pdo$ bundle $\calE$ 
is a $\Psi${\it DO connection} if its connection
one-form takes values in $\cl[](M,E)$ in any local trivialization.

\begin{remark}  If $\theta$ is the locally defined connection one-form
of a $\pdo$ connection on a bundle $\calE$ modeled on $C^\infty(M,E)$, 
then under a gauge change $g$, $\theta$
transforms to $g^{-1}\theta g + g^{-1}dg.$  Since $g^{-1}dg$ is zeroth
order if $g$ is nonconstant, the connection one-form is usually of
non-negative order.  (For left invariant connections on loop groups,
$g^{-1}dg$ vanishes, and $\theta$ can be of any order.)  When $\theta$
is non-positive order, it is a bounded operator and hence extends to a
connection on the extension of $\calE$ to an $H^s(M,E)$-bundle.  We
call a connection with connection one-form taking values in
$\cl[\leq 0](M,E)$ a {\it  $\cl[\leq 0]$-connection}.  The curvature
of such a connection is a $\cl[\leq 0]$-valued two-form on the base.
\end{remark}

\section{Traces, weighted traces, and corresponding  Chern
 classes in infinite dimensions}

In this section we examine two examples of Proposition 2.2.
The trace is furnished by the Wodzicki residue in the first example,
and by various traces applied to
the leading order symbol of a zeroth order $\pdo$ in the
second.  We will focus on the analytical/geometric aspects of this
approach, and leave topological aspects to a forthcoming paper.

In each case, we begin with a $\pdo$ bundle
$\calE$ over a base space $B$, with fiber modeled on $C^\infty(M,E)$
or $H^s(M, E)$, equipped with a $\pdo$ connection.
(We will let the connection one-form take values in  $\cl[](M,E)$ in \S3.3).
  In the language of \S2, we
pass from $\calE$ to the corresponding principal bundle $P =
P^{\calE}$ with fiber modeled on $\cl[0]^*$.
  $\AdP$ is then the
subbundle of $\End(\calE)$ consisting of endomorphisms in $\cl[\leq 0].$
In particular, the curvature $\Omega$ of the connection on $\calE$ takes values
in $\AdP$, and we can apply the Chern-Weil machinery of \S2 to
$\Omega.$  

Note that in this section, we are treating the structure group $\cl[0]^*$ as a
generalization of $GL(n,\C).$  As in finite dimensions, we focus only
on invariant polynomials on the Lie algebras given by traces.  We do
not discuss the interesting question of whether all such polynomials
on $\cl[\leq 0]$ are generated by these traces.

Exactly how these examples generalize the finite dimensional situation 
is open to interpretation.  When the manifold is reduced to a point,
the leading symbol of an endomorphism in the fiber, a ``zeroth order
$\pdo$,'' is just the endomorphism itself, and the only trace, up to
normalization, is the ordinary trace on a vector space.  In contrast,
in finite dimensions the Wodzicki residue vanishes.  So in this
interpretation, the Wodzicki residue is a purely infinite dimensional
phenomenon, while the symbol trace generalizes the finite dimensional
theory (since the leading symbol of a finite dimensional linear
transformation $A$ is just $A$).

On the other hand, both the Wodzicki residue and the symbol trace
appear in the most divergent term in asymptotic expansions: the
Wodzicki residue of an operator $A$ is the residue of the pole of the
zeta function regularization
$\Tr(AQ^{-s})$ at $s=0$ (for any positive elliptic operator
$Q$), and the symbol trace is related to the coefficient of the most
divergent term in the heat operator regularization
$\Tr(Ae^{-\e Q})$ as $\e\to 0.$ (The last statement is proved in
Proposition 3.4.)  Since the two corresponding
``regularizations'' in finite dimensions using a positive definite
matrix $Q$ simply reduce to $\Tr(A)$, we can alternatively view
both examples as proper generalizations of the finite dimensional
Chern-Weil theory.

There is a third approach to generalizing from finite dimensions,
namely by taking the finite part of a regularized trace.  In the zeta
and heat regularizations above, this amounts to picking a different
term in the Laurent/asymptotic expansion, one that is more difficult
to calculate in general.  This approach was taken in \cite{PR} with
limited success.  In \S3.3, we use weighted traces to produce closed
forms from the most divergent term in an asymptotic series.  This is
not an example of Proposition 2.2, as the weighted trace does not satisfy
(2.1).

In summary, not only is there no canonical generalization of finite
dimensional Chern-Weil theory, there is no canonical
interpretation of whether a specific method is indeed a proper
generalization.  The particular choice of regularization depends on a
combination of physical motivation, computability, and perhaps aesthetics.

\subsection{The Wodzicki residue}\label{wodzicki}

Recall that the Wodzicki residue $\res[A]$ 
of a $\pdo$ $A$ acting on sections of
a bundle $E$ over a closed manifold $M$ is defined to be the residue
of the pole term of $\Tr(AQ^{-s})$ at
$s=0$, for an elliptic operator $Q$ with certain technical conditions.
 Alternatively,
$\res[A]$ is proportional to  the coefficient of $\log \e$ in the asymptotic
expansion of $\Tr(Ae^{-\e Q})$ as $\e\to 0,$ for 
$Q\in \Ell^+(M,E).$  The strengths of the
Wodzicki residue are (i) its local nature:
\[ \res[A] = \frac{1}{(2\pi)^n}\int_{\sm}\tr\ \sigma^A_{-n}(x,\xi)\  d\xi
\ dx,\] 
where $n = {\rm dim}(M),$ $\sm$ is the unit cosphere bundle of $M$, and
$\sigma^A_{-n}$ is the $(-n)^{\rm th}$ homogeneous piece of the symbol
of $A$; and (ii) the fact that is it the unique trace on 
$\cl[] = \cl[](M,E)$,
up to normalization.  Its drawback is its
vanishing on all differential and multiplication
operators, all trace class operators
(and so all $\pdo$s of order less than $-n$) and all
operators of non-integral order.

Given an infinite dimensional bundle ${\mathcal E}$ over a base $B$
with fibers modeled either on $H^s(M,E)$
$E$ (with $s \gg 0$) or on $C^\infty(M,E)$, 
and a connection on ${\mathcal E}$ with
curvature $\Omega\in \Omega^2(B, \cl[])$, we can form the $k^{\rm th}$
Wodzicki-Chern form by setting
$$c_k^{\rm w}(\Omega) = \res[\Omega^k] \in \Omega^{2k}(B).$$
By Proposition 2.2, 
 $c_k^{\rm w}(\Omega)$ is closed and independent
of the connection. 

By the remark after Proposition 2.2, 
the Wodzicki-Chern classes vanish for
current groups. For completeness, we give a direct proof; see
\cite{jp} for a more general discussion.
For simplicity, we only consider Hilbert current groups,  $\calC = H^s(M, G)$,
 the space  of
$H^s$-maps from a closed Riemannian manifold $M$ to a semi-simple Lie group
$G$ of compact type.
(These assumptions ensure 
that the Killing form is nondegenerate and that the adjoint
representation 
is antisymmetric for this form).
$\calC$  is a Hilbert Lie group with     Lie algebra
$H^s(M, A)$, the space of $H^s$ sections
of the trivial bundle $M\times A$, where
$A = {\rm Lie}(G)$.
Thus the tangent bundle $T\calC$ is a $\pdo$ bundle with
fibers modeled on $H^s(M, A)$.
For $\Delta$ the Laplacian on functions on $M$, we set
 $Q_0:=\Delta \otimes 1_{A}, $
a second order elliptic operator acting
densely on $ H^s(M, A)$.
$Q_0$ is
non-negative for the scalar product
$\langle \cdot, \cdot \rangle_0:=
\int_M \dvol(x) (\cdot, \cdot)$,
where $(\cdot, \cdot)$ is minus
 the Killing form.
$T\calC$ has a left-invariant
weight $Q_\gamma= L_\gamma Q_0 L_\gamma^{-1}$
(i.e.~a family of elliptic operators on the fibers), where $L_\gamma$ is
left translation by $\gamma\in \calC$.

$\calC$ has  a  left-invariant
Sobolev $s$-metric   defined 
by
 $$\langle \cdot, \cdot \rangle^s:=
\langle {Q_0}^{s\over 2} \cdot,
{Q_0}^{s\over 2}\cdot
\rangle_0, $$
where ${Q_0}$ is really $ Q_0 + P$, for $P$ the orthogonal projection
of $Q_0$ onto its kernel.
The corresponding left-invariant Levi-Civita connection 
has the global expression
$\nabla^s   = d+ \theta^s$,
with $\theta^s$ a left-invariant
 $\End(T{\mathcal C})$-valued one-form on $\calC$  induced by 
the $\End(H^s(M,A))$-valued one-form on $H^s(M,A)$
\begin{equation}\label{lg}\theta_0^s(U)= {1\over 2} \left(\ad_U+
{Q_0}^{-s} \ad_U\
{Q_0}^s -
{Q_0}^{-s} \ad_{
{Q_0}^sU}\right),\end{equation}
 for $ U\in H^s(M,A)$ 
(see \cite[(1.9)]{F} up to a sign error).  $\theta^s$ takes values in
$\cl[\leq 0](M,M\times A).$  

There is another natural left-invariant connection on $\calC$ given by
\begin{equation} \label{another}\widetilde \nabla^s:= d+ \widetilde 
\theta^s,\end{equation}
defined by the corresponding
$\widetilde\theta_0^s(U) :=   {Q_0}^{- s} \ad_U
\  {Q_0}^{ s}.$   $\widetilde\theta^s$ also takes values in
$\cl[\leq 0](M,M\times A)$, and so $\nabla^s$ and $\widetilde\nabla^s$ can
be joined by a line of such connections.

We now show that the Wodzicki-Chern forms vanish for  
$\widetilde \nabla^s$,
from which the vanishing of the Wodzicki-Chern
classes of $\nabla^s$ follows.
The connection $\widetilde \nabla^s$ has curvature 
\begin{eqnarray*} \widetilde \Omega^s(U,V)&=& [\widetilde\theta_0^s(U),
\widetilde\theta_0^s(V)]- \widetilde\theta_0^s(
[U, V])
=  {Q_0}^{- s} [\ad_U, \ad_V]
  {Q_0}^{ s}- {Q_0}^{- s}  \ad_{[U, V]}
  {Q_0}^{ s}\\
&=&
{Q_0}^{- s}\left( \Omega^0(U, V) \right){Q_0}^{s},\end{eqnarray*}
for $ U,V\in
H^s(M, A)$, and where
$\Omega^0(U, V)=  [\ad_U, \ad_V]-
     \ad_{[U, V]}$ is a multiplication operator. As before, 
applying the Wodzicki residue to any analytic function $f$ yields
$$
\ress(f( \widetilde \Omega^s ))(U, V)= \ress \left( {Q_0}^{-
s}f\left(\widetilde
\Omega^0(U, V)
\right){Q_0}^{ s}\right)
= \ress \left(  f\left(\Omega^0(U,
V)
\right) \right)
=0,$$ since $\Omega^0(U,
V) $ is a multiplication operator.

In general,
no examples seem to be known where the Wodzicki-Chern classes 
are non-zero; see \cite{jp} for examples on spaces of maps
where the forms vanish.

\subsection{Leading symbol traces}

The uniqueness of the trace on $\cl[]$ defined by the Wodzicki residue
does not rule out the existence of other traces on subalgebras of
$\cl[]$.  Indeed, 
the ordinary operator trace on $\cl[< -n]$ is an example.
In this subsection, we will introduce a family of traces on $\cl[\leq 0]$
and show that they produce non-vanishing Chern classes on the
universal bundle associated to the gauge group for the $\calE = TLM$, the
tangent bundle
to the free loop space of a Riemannian manifold $M$. 
 To our knowledge, this is the first example of
non-vanishing Chern classes of infinite dimensional bundles above
$c_1.$  



We first produce a ``trace'' on $\cl[\leq 0]$ with values in $\sm$, and an
associated family of true traces.  Let $\calD(X)$ denoted the space of
complex valued distributions on a compact topological space $X$.

\begin{lemma}\label{symboltracelemma}
Let $\sigma_0^A$ be the zeroth order symbol of $A\in \cl[\leq 0].$ 
  The map  
$\Tr = \Tr_0: \cl[\leq 0](M,E) \to C^\infty (S^*M)$\ defined by\hfill \\
 $\Tr(A) = 
\tr_x( \sigma_0^A(x,\xi))$ has $\Tr(A+B) = \Tr(A) + \Tr(B)$, $\Tr(\lambda A )
= \lambda \Tr(A),$  and $\Tr(AB) =
\Tr(BA).$ For any $f\in \calD(\sm)$, the map $\Tr^f:\cl[\leq 0]\to\C$ given
by $\Tr^f(A) = f(\Tr(A))$ is a trace.
\end{lemma}

\begin{proof}  Certainly taking the $0^{\rm th}$ order symbol
is linear.  
Since $\sigma_0$ is multiplicative on $\cl[\leq 0]$, we have
$$\tr_x\
\sigma_{0}^{AB} = \tr_x(\sigma_{0}^A\cdot\sigma_{0}^B) = 
\tr_x(\sigma_{0}^B\cdot\sigma_{0}^A) = \tr_x\ \sigma_{0}^{BA}.$$  The
proof of the second statement is immediate.
\end{proof}

\begin{remarks}  (1)  There are corresponding traces on $\cl[\leq p]$ for
any $p<0.$ In fact, for $r\in [2p,p]$, $\Tr_r A=
\tr_x(\sigma_r^A(x,\xi))$ is a trace, as $\Tr_r(AB)$ trivially
vanishes for $r>2p.$  The proof of the lemma covers the case $r = 2p$.

(2) Let $Q\in\Ell^+(M,E)$ have scalar leading symbol
$\sigma_L^Q(x,\xi) = f(x,\xi){\rm Id}$.
Define $f\in\calD(\sm)$ by 
$$f(\phi) = \frac{\Gamma\left(\frac{n}{q}\right){\rm dim}(E)}{q^2 
(2\pi)^n}\int_{\sm}
\phi\cdot f(x,\xi)^{-{n\over q}}, $$
where $n={\rm dim}(M), q = {\rm ord}(Q).$  Then $\Tr^f(A)$ is the
leading term in the asymptotics of $\Tr(Ae^{-\e Q})$ (see Proposition
\ref{oldb}).

   \end{remarks}

Recall that the ring of characteristic classes for e.g. $U(n)$ bundles
is generated by the Chern classes $c_k = [\Tr(\Lambda^k \Omega)]$, or
equivalently by the components $\nu_k = [\Tr(\Omega^k)]$  of the
Chern character.  Note we are momentarily
distinguishing between $\Tr(\Lambda^k
A)$ and $\Tr(A^k)$ for a matrix $A$. We will concentrate on Chern
forms, and abuse notation by writing $c_k = [\Tr(\Omega^k)].$

\begin{definition}  Let $\calE$ be a bundle over $B$ modeled on
$H^s(M,E)$ or $C^\infty(M,E),$ and let
$\nabla$ be a connection on
$\calE$ with curvature $\Omega$ in $\cl[\leq 0](M,E).$
The $k^{\rm th}$ 
Chern class of $\nabla$
with respect to  $f\in
\calD(S^*M)$ is defined to be the de Rham cohomology class
\begin{equation}\label{def} [c_{k}^{f}(\Omega)] =
[f(\tr_x\sigma_0^{\Omega^k}(x,\xi))]\in H^{2k}(B;\C).\end{equation}
 \end{definition}

\begin{remarks} (1)  As an example, if $f\equiv 1\in C^\infty(\sm)$, then 
\[ c_k^1(\Omega) = \int_{\sm} \tr_x\ \sigma_0^{\Omega^k}(x,\xi).\]
At another extreme, if $f = \delta_{(x_0,\xi_0)}$ is a delta function,
then 
\[ c_k^f(\Omega) = \tr_x\ \sigma_0^{\Omega^k}(x_0,\xi_0).\]

(2) As in the previous remark, we can define Chern classes for connections
with curvature forms taking values in $\cl[\leq p]$ for any $r\in
[2p,p]$, $p<0.$
Note that for $r<p$ and e.g.~$f\equiv 1$, 
these classes are defined only after a
choice of coordinates on $M,E$ and a partition of unity on $M$, since
integrals of
non-leading order symbols depend on such choices.
\end{remarks}

Proposition 2.2
immediately gives the following:

\bigskip
\begin{theorem}  Let $ \nabla $ be a $\cl[\leq 0]$ connection on a $\pdo$
bundle $\calE$.
Then
the differential forms
$c_{k}^{f}(\Omega)$ 
are closed, and their cohomology classes are 
independent of  the choice of $\cl[\leq 0]$ connection. \end{theorem}

\begin{proof} By Lemma 3.1, 
$\Tr^f$ is a trace on $\cl[\leq 0](M, E)$, so we can apply
 Lemma 2.1 to the principal bundle $P= P^{\cal E}$ built from
 ${\cal E}$ to get the relation $d\circ \Tr^f=\Tr^f\circ \nabla^\ad$.
We then apply Proposition 2.2 to get the corresponding Chern
 classes $[c_k^f(\Omega)]$.  \end{proof}

\begin{remark} For the  linear functionals $\Tr_p$ on $\cl[\leq p]$, the proof 
that the Chern forms are closed goes through, but the proof of their
independence of choice of connection breaks down.  The point is that
the line of $\cl[\leq 0]$ 
connections joining two connections with curvature forms
lying in $\cl[\leq p]$ might not stay in this class. \end{remark}

When the structure group reduces to a gauge
group, we can construct an example of non-zero Chern classes
$[c_k^{f}(\Omega)].$ Fix $n>k$ and consider the Grassmanian $BU(n) =
G(n,\infty)$ with its universal vector bundle $\gamma_n$. We consider
the pullback bundle $E =\pi^*\gamma_n$ over $S^1\times BU(n)$, with
$\pi$ the projection onto $BU(n).$ We now ``loopify'' to form $B =
L(S^1\times BU(n))$, the free loop space of $S^1\times BU(n)$, with
bundle $\calE$ whose fiber over a loop $\gamma$ is the space of smooth
sections of $\gamma^*E$ over $S^1$.  ($\calE_\gamma$ is the space of
loops in $E$ lying over $\gamma$, suitably interpreted at
self-intersection points of $\gamma$, so we will write $\calE =
L\pi^*\gamma_n$.)  Since $\gamma^*E$ is (non-canonically) isomorphic
to the trivial bundle $S^1\times \BbC^n$ over $S^1$, it is easily
checked that the structure group for $\calE$ is the gauge group
$\calG$ of this trivial bundle.

 Take a hermitian
connection $\nabla$ on $\gamma_n$ (e.g.~the universal
connection $PdP$, where $P_x$ is the projection of $\BbC^\infty$ onto
the $n$-plane $x$)  and its pullback connection 
$\pi^*\nabla$ on $S^1\times BU(n).$
As in the case of the tangent
bundle to a loop space, we can take an
$L^2$ or pointwise connection $\nabla^0$ on $\calE$ by setting
$$\nabla^0_XY(\gamma)(\theta) = 
\pi^*\nabla_{X(\theta)}Y(\theta),$$
for $X$ a vector field along $\gamma$ (i.e. a tangent vector in $B$ at
$\gamma$) and $Y$ a local section of $\calE.$  
The curvature
$\Omega^0$ acts pointwise and hence is a multiplication operator.
In particular, its symbol is independent
of $\xi$. 

Pick the distribution $\delta = (1,+)$ on
$C^\infty(S^*S^1) = C^\infty(S^1\times \{\pm\partial_\theta\})$: i.e.
$$\delta(f(\theta,\partial_\theta), g(\theta,-\partial_\theta)) =
\frac{1}{2\pi}\int_{S^1} f(\theta,\partial_\theta)\ d\theta.$$ 

We claim that $[c_{k}^{\delta}(\Omega^0)]$ is nonzero in
$H^{2k}(B; \C).$  To see this, let $a = a_{2k}\in H_{2k}(BU(n),\C)$ 
be such that
$\langle c_k(\gamma_n),a \rangle = 1.$  
Define $c\in H^{2k}(B)$ to be $c = \beta_*a$,
where $\beta:BU(n) \to L(S^1\times BU(n))$ is
given by $\beta(x)(\theta) = (\theta,x)$.  
Now 
\begin{equation}\label{nine}
\langle  [c_k^{\delta}(\Omega^0)], c\rangle =  
\langle  [c_k^{\delta}(\Omega^0)], \beta_*a \rangle =  
\langle [\beta^* c_k^{\delta}(\Omega^0)], a \rangle,\end{equation}
since $\beta$ has degree one. 
For $\gamma\in L(S^1\times BU(n))$, we have
   \begin{equation}\label{ten}
   c_k^{\delta}(\Omega^0)(\gamma) =\frac{1}{2\pi}\int_{S^1}
 \tr\left(\sigma_0^{(\Omega^0)^k} (\gamma(\theta),
\partial_{\gamma(\theta)}) \right)\ d\theta
= \frac{1}{2\pi}\int_{S^1}
 \tr\left(\Omega^k_{\gamma(\theta)}\right)\ d\theta.\end{equation}
For a tangent vector $X\in T_xBU(n)$, it is immediate that
$\beta_*(X)\in  T_{\gamma_x}B$ has 
$\beta_*(X)(\theta,x) = (0,X).$  Thus by (\ref{nine}),
\begin{equation}\label{eleven}
\beta^* c_k^{\delta}(\Omega^0) (X_1,\ldots,X_{2k}) =
\frac{1}{2\pi}\int_{S^1} 
\tr(\Omega^k)(X_1,\ldots,X_{2k})
= \tr(\Omega^k)(X_1,\ldots,X_{2k}).\end{equation}
Combining (\ref{ten}) and (\ref{eleven}), we get
$$\langle  [c_k^{\delta}(\Omega^0)], c\rangle =  \langle
[\tr(\Omega^k)], a \rangle = 1.$$

In particular, the class $[c_k^{\delta}(\Omega^s)]$ must be
non-zero.  

\begin{theorem}  \label{nonzero}
The cohomology classes $ [c_k^{\delta}(\Omega)] $
are
non-zero in general.  In particular, the corresponding classes for
the universal  bundle $E\calG$ are nonzero 
in the
cohomology of the classifying space
$B\calG$, where $\calG$ is the gauge group of the
trivial bundle $S^1\times \C^n$ over $S^1.$ \end{theorem}
\bigskip

We have shown the first statement.
To explain the second statement, note that
although the structure group of $L\pi^*\gamma_n$ is the
gauge group of the trivial bundle $S^1\times \C^n$ over $S^1$, the
curvature of the connection will take values in $\cl[0](S^1,S^1
\times\C^n)$ of this
bundle.  As a result, the classifying space is really 
$B{\mathcal Cl}^*_{0}$.
However, the principal symbol map gives a retraction of ${\mathcal
 Cl}^*_{0}$  onto the gauge group of the trivial bundle over
$S^*S^1$, which is just two
 copies of $\calG.$  Thus $B{\mathcal Cl}^*_{0}$ is homotopy
equivalent
to $ B\calG \coprod
 B\calG$, and each $[c_k^{(1,\pm)}]$ is non-zero in one copy of $B\calG.$  
 The proof of the second statement depends on the existence of a
universal connection on $E\calG$ over $B\calG$, and will be given
in \S4.

In fact, $B\calG$ equals $L_0BU(n)$,
the space of contractible loops on $BU(n)$ \cite{AB}.  It is known that
 $H^*(B\calG,\C)$ is a super-polynomial 
(i.e.~super-commutative) algebra
 with one generator in each degree $k\in \{1,\ldots,2n\}.$  
In analogy with finite dimensions, we
 conjecture that
 $[c_k^\delta]$ is a nonzero multiple of the generator
in degree $2k.$  For $k=1$, this is clear.

We now outline a conjectured construction of geometric representatives
of the odd generators in $H^*(B\calG).$  The tangent bundle $TLM$ of
any loop space splits off a trivial line bundle, namely the span of
$\dot\gamma$ at the loop $\gamma.$   Note that for any connection on a
bundle over $LM$, we have
\begin{equation}\label{lie}
di_{\dot\gamma}c_k^f(\Omega)  =
di_{\dot\gamma}c_k^f(\Omega) + i_{\dot\gamma}dc_k^f(\Omega) = 
 L_{\dot\gamma}c_k^f(\Omega),\end{equation}
where $i$ is interior product and $L$ is Lie derivative.  We state
without (the elementary) proof that (\ref{lie}) implies
\begin{equation}\label{lie2} 
di_{\dot\gamma}c_k^f(\Omega)  = c_k^{\partial
f}(\Omega),\end{equation}
where $\partial f$ is the derivative of $f$ as a distribution on
$S^*S^1.$  In particular, we see that
$i_{\dot\gamma} c_k^{(1,\pm)}(\Omega)$ are closed forms.
It remains to be seen if $[i_{\dot\gamma}
c_k^{(1,\pm)}(\Omega)] \in H^{2k-1}(B\calG \coprod B\calG)$ are non-zero.

We can also use (\ref{lie2}) to understand the dependence of
$[c_k^f(\Omega)]$ on $f$.  Since $f$ is a zero current on $S^*S^1$ and
hence is trivially closed, and since exact zero currents $\partial f$
produce vanishing Chern classes by (\ref{lie2}),
we see that the
space of classes $$\{[c_k^f(\Omega)]:f\in \calD (S^*S^1)\}\in
H^{2k}(B\calG)$$ is isomorphic to the zeroth cohomology group of
complex currents on $S^*S^1.$ 
(Here we are extending the usual 
confusion of functions $f$ and one-forms $fd\theta$ on $S^1$ to a
confusion of zero- and one-currents.)
This cohomology group is isomorphic to $H_0(S^*S^1)$, and is spanned
by $(1,\pm).$ (The reader may wish to check directly from (\ref{lie2})
that all $\delta$-function currents on one copy of $S^1$ produce the
same cohomology class.  A more interesting exercise is to show that
these delta functions produce the same cohomology class as one of
$(1,\pm)$ using Fourier series.)

\begin{remarks}
(1)
The general  case, where the
gauge group is associated to the bundle $E$ over a closed manifold $M$, is more
complicated.  
The cohomology
of $B\calG$ is known, and in general has odd dimensional
cohomology \cite{AB}.  We do not know at present
 which part of $H^*(B\calG)$ is
spanned by $[c_k^f(\Omega^{E\calG})]$, where we use the universal
connection mentioned above.  We also do not know how to produce
geometric representatives of odd dimensional classes in $H^*(B\calG)$,
nor do we know how the Chern classes depend on the distribution.

(2)
In the loop group case, Freed showed \cite{F} that the curvature
$\Omega^s$ of
the $H^s$ Levi-Civita connection is a $\pdo$ of order $-1$ for $s>
1/2.$  It is for this reason that we did not use this Levi-Civita
connection in our construction; the  Chern forms built from $\sigma_0$
trivially vanish.
\end{remarks}

\subsection{Weighted first Chern forms on loop groups}

In this subsection we treat the weighted first Chern forms of
\cite{PR} for loop groups.  
In contrast to the traces in previous sections, 
the weighted traces are associated to the finite part of an
asymptotic expansion, and do
not produce closed characteristic forms
in general.  There is an important exception:
 for the $H^{1/2}$ metric on a loop group, 
Freed \cite{F} showed that a conditional trace (recalled below) of the
curvature form converges and equals the K\"ahler form.  Thus this
conditional first Chern form
 is certainly closed with non-zero cohomology class.  In
\cite{CDMP}, it was shown that this curvature form is a weighted first
Chern form.

We first give an alternative proof that the conditional first Chern
form is
closed (Proposition \ref{last2}).  We then show 
more generally that forms occuring as the most
divergent term in an asymptotic expansion associated to 
weighted traces are closed under certain hypotheses (Theorem
\ref{last3}). In particular, we recognize the most divergent term in
Proposition \ref{oldb} as a symbol trace from \S3.2, which we know
leads to closed characteristic forms.
The conditional first Chern form
case is special precisely because the most divergent term and
the finite part coincide on $LG.$  Thus, Freed's computation 
fits into
the framework of weighted traces, which
is strengthened in \cite{PR2}.

We first recall the definition of  weighted traces. These
are  functionals on $\cl[](M,E)$ defined for
$Q\in\Ell^+(M,E)$ by $\tr_{\e}^Q(A)
:= \Tr(Ae^{-\e Q})$ for any $\e>0.$  Recall that a $\pdo$ bundle $\calE$ 
with structure group 
$\cl[0]^*(M,E)$ has an associated bundle of
algebras $\cl[\leq 0](\calE) = \AdP^{\calE}$ 
with fibers modeled on 
$\cl[\leq 0](M,E).$   
A weight is a section $Q\in\Gamma(\calB)$ with $Q$ elliptic with
positive definite leading symbol and of constant order.  These
conditions are independent of local chart, since the transition maps
are $\pdo$s.
In particular,
if $\{g_b\}$ is the transition map between two trivializations of
$\calE$ over $b$, then $Q_b$ transforms into $g_b^{-1}Q_b g_b;$ the
same holds for sections of $\calB$.  For $A\in\Gamma(\calB)$,
$\tr_\e^Q(A)$ is well-defined, since 
\begin{equation}\label{trq}
\tr_\e^{g^{-1}Qg}(g^{-1}Ag) = \Tr(g^{-1}Age^{-\e g^{-1}Qg}) =
\Tr(g^{-1}Agg^{-1} e^{-\e Q}g) = \tr_\e^Q(A).\end{equation} 
We set $\tr^Q(A)$ to be the finite part of $\tr_\e^Q(A)$ as $\e\to 0.$
In other words, $\tr^Q(A)$ is the coefficient of $\e^0$ in the
asymptotic expansion (\ref{bone}).  This is equivalent to taking the
zeta function regularization $\Tr(AQ^{-z})|_{z=0}$, provided $Q$ is
invertible and (\ref{bone}) contains no log terms.

As a preliminary result,
we show that the most divergent term in the
asymptotics of $\tr(Ae^{-\e Q})$ is given by an integral of the
leading order symbol $\sigma_L(A)$ of $A$ if $Q$ has leading order symbol 
$\sigma_L(Q)(x,\xi) = ||\xi||^k$ for some $k$. 
 Here $A$ is a zeroth order
$\pdo$ and $Q$ is a Laplacian-type operator. 
This folklore result follows from the analysis developed in \cite{GS},
while the analysis in the following proof is hidden in the local
nature of the Wodzicki residue.

\begin{proposition} \label{oldb}
Let  $A\in \cl[\leq 0](M, E)$ have integral order $a \geq -n, n = {\rm
dim}(M)$, and let  $Q$ be
 an elliptic $\pdo$ of order $q$ with positive scalar leading symbol
$\sigma_L(Q)(x,\xi) = f(x,\xi){\rm Id}.$ Let $c = c(n,a,q)$ be 
$$c = { \Gamma({n+a\over q}){\rm dim}(E)(n-1)! \over
q(2\pi)^n}, \ {\rm if}\ a\neq -n;\ \  c = \frac{{\rm
dim}(E)(n-1)!}{(2\pi)^n}, \ {\rm if}\ a = -n.$$
Then as $\e\to 0$,
$$\tr(A e^{-\e Q}) =  
c\int_{S^*M} \tr \left(\sigma_a(A) \right)
\left(f(x,\xi)\right)^{-\frac{n+a}{q}}
 \cdot \e^{-{n+a\over q}} + \oo(\e^{-{n+a\over q}}).$$ 
In particular, if
$\sigma_L(Q)(x,\xi) = ||\xi||^k$ for some $k$, then
$$\tr(A e^{-\e Q}) =  
c\int_{S^*M} \tr \left(\sigma_a(A)\right) 
 \cdot \e^{-\frac{n+a}{q}}+ \oo(\e^{-{n+a\over q}}).$$ 
\end{proposition}

\begin{proof} We want to compute the coefficient  $a_0(A, Q)$
in the known asymptotic expansion 
\begin{equation}\label{bone}\tr(A e^{-\e Q})=
\sum_{j=0}^{a+n} a_j(A, Q) \e^{{j-a-n\over q}}+ b_0(A, Q) \log \e
+\OO(1),\end{equation}
for general $A\in\cl[](M,E)$ of order $a$, and with
$a_j(A, Q), b_0 (A, Q)\in \C$. The coefficient $b_0(A,
Q)$ satisfies $b_0(A, Q)= -{1\over
q}\ress(A)$.

 A  Mellin transform yields 
$$a_0(A, Q)= 
 {\rm Res}_{z={ n+a\over q}} \Gamma\left({n+a\over q}\right)\tr(A
Q^{-z}),$$ 
(the case $A=1$ considered in \cite[(12)]{K} easily
extends to a general $\pdo$ $A$). 
Thus, for $A$ as in the hypothesis, 
$$a_0(A, Q)= {\Gamma({n+a\over q}) \over q\Gamma(n+a) } a_0(A, Q^{1\over
q}),\ {\rm if}\ a\neq -n;\ a_0(A, Q)= a_0(A, Q^{1\over q}), \ {\rm
if}\ a= -n.$$ 
Thus it suffices to prove the formula for 
$Q_1$ of order one. 
Since ${\rm ord}(AQ_1^{-(n+a)}) = -n$,  (\ref{bone}) becomes
 $$\tr(AQ_1^{-(n+a)} e^{-\e Q_1})=-\ress(AQ_1^{-(n+a)}) \log \e+\OO(1).$$
Differentiating this expansion
$n+a$ times  (valid since $n+a\geq 0$) with respect to $\e$, we get
$$\tr(A e^{-\e Q_1})\sim (n+a-1)!\ \ress(AQ_1^{-(n+a)})  \e^{-(n+a)}.$$
The local formula for the Wodzicki residue yields: 
\begin{eqnarray*} \ress(AQ_1^{-(n+a)}) &=&
 {1\over (2\pi)^n} \int_{S^*M} \tr(\sigma_{-n} (A
Q_1^{-(n+a)}))\\
&=& {1\over (2\pi)^n} \int_{S^*M} \tr\left(\sigma_{a}
(A)\sigma_{-(n+a)} (Q_1^{-(n+a)})\right)\\
&=& {{\rm dim}(E)\over (2\pi)^n} \int_{S^*M} \tr(\sigma_a (A))
f(x,\xi)^{-(n+a)}.\end{eqnarray*}
Hence our original $Q$ has
$$a_0(A, Q)= {\Gamma({n+a\over q}) \over q\Gamma(n+a) } a_0(A, Q^{1\over
q})= c \int_{S^*M} \tr\left(\sigma_a
(A)\right)\left(f(x,\xi)\right)^{-(n+a)/q},$$ 
if $a \neq -n$, and there is a similar formula for $a = -n.$
\end{proof}

We now show that the weighted first Chern forms 
are closed.  We use the notation of the previous section.

\begin{proposition}\label{last2} Let $\Omega = \Omega^{(1/2)}$ 
be the curvature of the Levi-Civita
connection $\nabla = d + \theta^{(1/2)}$ for the $H^{1/2}$ metric on the loop
group $LG$.
Then the weighted first Chern form $\tr^Q(\Omega)$ is
closed. \end{proposition} 

\begin{proof} The conditional trace of the curvature in \cite{F} is
$\tr(\tr_\lie (\Omega))$, where $\tr_\lie$ denotes the trace with
respect to the Killing form in the Lie alebra of $G$, and the outer
trace is the ordinary operator trace.  In particular, $\tr_\lie
(\Omega)$ is a trace class $\pdo$ on the trivial $\C^{{\rm dim}(G)}$
bundle over $S^1.$
It is shown in \cite{CDMP} that
$\tr^Q(\Omega) = \tr(\tr_\lie (\Omega)).$  Since $\tr_\lie
(\Omega)$ is trace class, we have
$$\tr(\tr_\lie (\Omega)) = \lim_{\e\to 0} 
\tr\left(\tr_\lie (\Omega)e^{-\e Q}\right) = \lim_{\e\to 0}
\tr(\Omega e^{-\e Q}).$$
Since $\Omega$ is a $\pdo$ of order $-1$ \cite[Thm.~1.11]{F}, the 
previous proposition gives
\begin{equation}\label{last4} \tr^Q(\Omega) = \lim_{\e\to 0} 
\tr\left(\Omega e^{-\e Q}\right) = a_0(\Omega,Q) =
\frac{{\rm
dim}(G)}{2\pi}\int_{\ssi} \tr_\lie
\sigma_{-1}(\Omega).\end{equation}
Thus 
$$d\tr^Q(\Omega) = \frac{{\rm
dim}(G)}{2\pi}\int_{\ssi} \tr_\lie \sigma_{-1}(d\Omega).$$
Since the connection one-form $\theta = \theta^{(1/2)}$ has order zero,
$\sigma_{-1}([\theta,\Omega])$ equals the leading symbol
$ \sigma_{L}([\theta,\Omega])$, which as usual vanishes.  Therefore
$$d\tr^Q(\Omega) = \frac{{\rm
dim}(G)}{2\pi}\int_{\ssi} \tr_\lie \sigma_{-1}(\nabla\Omega)=0.$$
\end{proof}

\begin{remark} It is not true that the K\"ahler class
$[\tr^Q(\Omega^{(1/2)})]$ is independent of the connection on $LG.$  Indeed,
replacing $\theta^{(1/2)}$ by $t\theta^{(1/2)},
 t\in [0,1]$, joins $\nabla$ to the
trivial connection $d$, which of course has $[\tr^Q(\Omega)] =0.$
The difficulty is that the order of the corresponding curvatures
$\Omega_t$ jumps from $-1$ at $t=1$ to zero for $t<1$.  This prevents
one from fixing a symbol term as in the proof of Proposition 3.5 for
all $t$. \end{remark}

We now generalize Proposition \ref{last2} to show that the most
divergent term in an asymptotic series associated to a weighted trace
is closed.

\begin{theorem}\label{last3}  Let $\calE\to B$ be a $\pdo$ bundle with
a $\cl[\leq 0]$ connection $\nabla$. 
  Let $Q\in\cl[](\calE)$ be a smooth family of
weights of constant order $q$ and 
with scalar leading symbol $\sigma_L(Q_b)$ independent of  $b\in
B.$  Let $A\in \Omega^k(B,\cl[](\calE))$ be a $\cl[](\calE)$-valued
form whose order $a\in \Z, a\geq -n, n = {\rm dim}(M)$ 
is  independent of  $b\in B$, and
such that $\nabla A =0.$  Write
$$\tr(A e^{-\e Q})=
\sum_{j=0}^{a+n} a_j(A, Q) \e^{{j-a-n\over q}}+ b_0(A, Q) \log \e
+\OO(1).$$
Then as elements of $\Omega^k(B,\C)$, 

(i) $\ress(A)$ is closed; 

(ii) $a_0(A,Q)$ is closed, if $a >-n$;

(iii) $a_0(A,Q)$ is closed, if $a= -n$ and $\ress(A) = 0$.
   \end{theorem} 
 
Note that the condition on the leading symbol is independent of
trivialization of $\calE.$  

\medskip

\begin{proof}  (i) This was shown in \S3.1.

(ii) and (iii)\  Under the hypotheses, $a_0(A,Q)$ is the most
divergent term in the asymptotic expansion.  For $c$ as in Proposition
\ref{oldb}, we have
\begin{eqnarray*} da_0(A,Q) &=&d\left[ c\int_{S^*M}
\tr\left(\sigma_a(A)\right) (f(x,\xi))^{-(n+a)/q}\right]\\
&=& c \int_{S^*M}\tr\left(d\sigma_a(A)\right) (f(x,\xi))^{-(n+a)/q}\\
&=& c \int_{S^*M}\tr\left(\sigma_a(dA)\right)
(f(x,\xi))^{-(n+a)/q}\\ 
&=& c \int_{S^*M}\tr\left(\sigma_a(dA + [\theta,A])\right)
(f(x,\xi))^{-(n+a)/q}\\  
&=& c \int_{S^*M}\tr\left(\sigma_a(\nabla A)\right)
(f(x,\xi))^{-(n+a)/q}\\  
&=&0.\end{eqnarray*}
Here we use the fact that $\nabla$ has non-positive order, so that 
$$\tr\left(\sigma_a([\theta,A])\right) = \tr\sigma_L([\theta,A]) =
\tr\left([\sigma_0(\theta),\sigma_a(A)]\right) = 0,$$
if $[\theta,A]$ has expected order $a$. Finally, $\sigma_a([\theta,A])=0$
trivially if the order of $[\theta,A]$ is less than $a$. \end{proof}

We now apply Theorem \ref{last3} to $LG$ to recover Proposition
\ref{last2}. Let 
$\Omega = \Omega^{(1/2)}, \theta = \theta^{(1/2)}.$
Note that
by (\ref{last4}), 
$\tr^Q(\Omega)$ is the most divergent term
$\alpha_0(\Omega,Q)$ 
(and the
finite part) of $\tr\left(\Omega e^{-\e Q}\right)$.
%
For $s = 1/2$, $a= -n$ 
and by case (iii), we must verify that $\ress(\Omega(X,Y)) =
0,$ for tangent
vectors $X,Y$ to $LG$.  Now 
$$\ress( \Omega(X,Y)) =
 \ress([\theta(X),\theta(Y)]) -
\ress(\theta[X,Y])
= -\ress(\theta[X,Y]).$$
Setting $U = [X,Y]$, we have by (\ref{lg})
\begin{eqnarray*} \ress(\theta(U)) &=&
 \frac{1}{2}\left( \ress(\ad_U)
+ \ress(Q^{-1/2}\ad_U Q^{1/2}) - \ress(Q^{-1/2}\ad_{Q^{1/2}U})\right)\\
&=& -\frac{1}{2} (\ress(Q^{-1/2}\ad_{Q^{1/2}U})),\end{eqnarray*}
since $\ress(\ad_U) =  \ress(Q^{-1/2}\ad_U Q^{1/2})$, and the Wodzicki
residue of this multiplication operator vanishes.  
Moreover, 
\begin{eqnarray*}\ress(Q^{-1/2}\ad_{Q^{1/2}U}) &=& \frac{1}{2\pi} \int_{S^*S^1}
\tr_\lie\left(\sigma_{-1}(Q^{-1/2} \ad_{Q^{1/2}U})\right)\\
& =& 
\frac{1}{2\pi} \int_{S^*S^1}\sigma_{-1}(Q^{-1/2})\tr_\lie(\ad_{Q^{1/2}U})\\
& =& 0.\end{eqnarray*}
Here we have used that $Q^{-1/2}$ has order $-1$, 
$\ad_{Q^{1/2}U}$ is a multiplication operator, and
$\tr_\lie(\ad_Z)=0$ for all $Z$ in the Lie algebra.

\begin{remark}  If the order of $A$ is zero, Proposition \ref{oldb}
implies that the most divergent term $a_0(A,Q)$ is a symbol trace in
the sense of \S3.2.  Thus Theorem 3.6 gives another proof that symbol
traces produce closed forms.  \end{remark}

\section{Universal connections}

In this section, we prove a Narasimhan-Ramanan theorem for universal
connections
on bundles  with structure group given by the gauge group $\calG$ of a
fixed Hermitian bundle $E$ over a closed manifold $M$.
 As explained
in the last section, this
justifies the second statement in Theorem 3.3.

We first show that the structure group of a $\pdo$ bundle 
reduces to a gauge group.
We have the  exact sequence at the (Lie)
algebra level 
$$0\to \cl[<0]\to \cl[\leq 0]\stackrel{\sigma_0}
{\to} C^\infty(S^*M,\End(\pi^*E))\to 0.$$
$\sigma_0$ is the zeroth order
symbol map; because zero order symbols have homogeneity zero in the
cotangent variable, these symbols are in fact functions on the unit
cosphere bundle
$S^*M$ with
values in $\End(E).$ This function space is the Lie algebra of the gauge
group $C^\infty(S^*M,\Aut(\pi^*E))$ 
of the pullback bundle $\pi^*E$ over $S^*M$, where $\pi:T^*M\to
M$ is the projection.
  We exponentiate 
this sequence 
to obtain 
\begin{equation}\label{split}
1\to 1 + \cl[<0] \to \Ell_{0}^*
\stackrel{\stackrel{\alpha}{\longleftarrow}}{\longrightarrow}
 C^\infty(S^*M,\Aut(\pi^*E)) \to
1,\end{equation} 
where $\alpha$ is a non-canonical splitting recalled below.
Since
$\cl[<0]$ is a vector space and hence contractible, and since it is
easy to verify that this exact sequence is a fibration, we have a 
 retraction of
$\Ell_{0}^* = \cl[0]^*$  onto 
$\alpha\left(C^\infty(S^*M,\Aut(\pi^*E))\right)$.

We now show that in fact
$\Ell_{0}^*$ deformation retracts onto this image.  This implies that
the corresponding classifying spaces satisfy $B\Ell_0^* \sim B
C^\infty(S^*M,\Aut(\pi^*E))$.  Thus we can take a gauge group as
structure group, although it is the gauge group for a bundle over
$S^*M.$

For the deformation retraction, it is standard that for $A\in
\cl[0]^*$, we can construct $\tilde B\in \cl[0]^*$ with the total
symbol of $\tilde B$ satisfying $\sigma_{\rm tot}(\tilde B) =
\sigma_0(\tilde B) = \sigma_0(A).$  In fact, if we fix a smooth
function $\psi:\R^n\to \R$ with $\psi\geq 0, \psi(0) = 0,$ and
$\psi(\xi) = 1$ for $|\xi| \geq 1,$ then we can set $\tilde B =
\alpha(A)$, where
$\alpha: C^\infty(S^*M,\Aut(\pi^*E))\to {\mathcal B}$ is defined by 
$$\alpha(\sigma)(f)(x) = \int_Me^{ix\cdot\xi}\psi(\xi)\sigma(x,\xi)\hat f(\xi)
d\xi.$$ 
Here, as usual, we have omitted notation for the trivialization of $E$
and partitions of unity.  
This is the splitting map $\alpha$ in
(\ref{split}). 

  Set $B = \tilde B + P_{\tilde B}$, where $P_{\tilde
B}$ is the orthogonal projection onto the kernel of $\tilde B.$  Note
that $\sigma_{\rm tot}(B) = \sigma_0(B) = \sigma_0(A).$  
$B$ is invertible, because it is an injective operator with ${\rm
index}(B) = {\rm index}(A) = 0.$  Thus
we can write $A = B(1+K)$, where $K$ is the
compact operator $B^{-1}(A-B).$  The set of $\lambda \in \C$ for which
$B(1 + \lambda K)$ is not injective is discrete, accumulates only at
infinity, and does not contain $0$
or $1$.  Since the same is true for $(1+\lambda K^*)B^*$, the union of
these two sets is a discrete set $S$ in $\C\setminus \{0,1\}$, and for
$\mu\in \C\setminus S$,
$B(1+\mu K)$ is invertible.  Fix a 
continuous curve $\gamma:[0,1]\to\C\setminus S, \gamma(0) = 0,
\gamma(1) = 1$, chosen so that $\gamma(t)=t$ except near
points of $S$, where it avoids points of $S$ via e.g.~a canonical
small semicircle in the upper half plane.  Then 
the family $A_t = B(1+\gamma(t)K)$ provides a 
deformation retraction of $\Ell_0^*$ to ${\mathcal B} = \{B\in
\Ell_0^*: \sigma_{\rm tot}(B) = \sigma_0(B)\}$.  
Since
elements of ${\mathcal B}$ are determined by a zeroth order symbol and
a smoothing operator, or equivalently by such a symbol and a smooth
kernel, we have
$${\mathcal B}  = \alpha(C^\infty(S^*M,\Aut(\pi^*E))) \times
C^\infty(M\times M, E\otimes E). $$
Since the set of smooth sections of a vector bundle
is contractible, we see that
$\Ell_0^*$ deformations retracts onto $ \alpha
( C^\infty(S^*M,\Aut(\pi^*E))).$

As a result, we may assume that our bundles have the gauge group
$  C^\infty(S^*M,\Aut(\pi^*E))$ as structure group.  We now begin the
demonstration  that the universal bundle for the structure group
$\calG = C^\infty(N,\Aut(F))$ has a universal connection.  Here the
automorphism group of a fiber is $G= U(n)$, but the argument applies
to the other classical groups.
The reader should
keep in mind our case
$S^*M = N, \pi^*E = F.$

Let $X$ be a (for us always locally trivial) principal $\calG$-bundle
over a topological space $B$.  
$X$ induces a $G$ principal bundle
$X'$ over $B\times N$ as follows.  We first set up a 1-1
correspondance of local sections $\Gamma_{loc}(X)\leftrightarrow
\Gamma_{loc}(X')$ by picking $U_i\subset B$ over which $X'|_{U_i}
\approx U_i\times \calG.$ For $s\in \Gamma(U_i)$, we map
$s\leftrightarrow u_s$, with $u_s(b,x) = s(b)(x).$ This agrees on
overlaps: on $U_i\cap U_j$, we have $s_i(b) = g_{ij}(b)s_j(b),
g_{ij}(b)\in \calG,$ and so
$$u_{s,i}(b,x) = g_{ij}(b,x)u_{s,j}(b,x),$$
with $g_{ij}(b,x) \equiv g_{ij}(b)(x).$ Now we define $X'$ over
$B\times N$ to be the $G$ bundle with transition functions
$g_{ij}(b,x)$ over $(U_i\cap U_j)\times N.$ In particular, we may
associate a $G$ bundle $E\calG'$ over $B\calG\times N$ to the universal
bundle $E\calG\to B\calG.$  Note that $E\calG = (\pi_1)_\ast E\calG'$
for the projection $\pi_1:B\calG\times N \to B\calG.$

We will use the following diagram to move from the universal
connection on $EG$ to a universal connection on $E\calG$.  
$$\begin{array}{ccccc}E\calG\simeq
(\pi_1)_\ast E\calG'& &E\calG'\stackrel{\alpha}{\simeq}
 \ev^*EG&\rightarrow&EG\\
\downarrow&&\downarrow&&\downarrow\\
B\calG&\stackrel{\pi_1}{\leftarrow}& B\calG\times
N&\stackrel{\ev}{\rightarrow}&BG\end{array}$$ 
The middle
isomorphism is contained in the proof of Lemma 4.1.

As a preliminary step,
we show how to classify connections on
$G$ bundles over $B\times N.$ 
 Recall that $B\calG = \MMaps_F(N,BG)$, the space
of maps which induce $F$ \cite{AB}.
  This gives the evaluation map $B\calG\times
N\to BG,\ \ev(f,x) = f(x).$ Note that $\ev^*EG$ restricted to
$\{f\}\times N$ is precisely $f^*EG$, which is isomorphic to the
unitary frame bundle $\Fr(F).$

\begin{remark}  By \cite{AB},
$E\calG =
\MMaps^{G}(\Fr(F), EG),$
 with the fiber over $f$ consisting of
bundle maps lying over $f$.  Here the superscript $G$ indicates
that the maps   are $G$-equivariant.
As a result, $E\calG$ consists of all
isomorphisms of $\Fr(F)$ with pullback bundles $f^*EG$.  Since every
connection on $\Fr(F)$ is isomorphic to the pullback $f^*\nabla^u$ of
the universal connection $\nabla^u$ on $EG$ for some $f\in B\calG$
\cite{NR},
we see that $E\calG$ contains all connections on $\Fr(F)$.  Connections
occur with repetitions if $\calG$ does not act freely on $\calAA$, the
space of all connections on $F$. (This
is the difference between using $\calG$ and the set $\calG_0$ of gauge
transformations which are the identity at a fixed point; $\calG_0$
acts freely on $\calAA$,
and so $E\calG_0 = \calAA.$)
\end{remark}

Let $(X',\nabla')$ be a $G$ bundle with connection over
$B\times N$, and let
$\beta:B\times N\to BG$ be a geometric classifying map --
i.e.~$\beta^*\nabla^u = \nabla'$, where we suppress the isomorphism
between $X'$ and $\beta^*EG.$ Set $$\beta_1:B\to B\calG =
\MMaps_F(M,BG), \ \ \beta_1(b)(x) = \beta(b,x).$$ Thus
$(\beta_1,\id):B\times N\to B\calG\times N.$ Since $\ev\circ(\beta_1,\id)
= \beta$, we have $$X' = \beta^*EG = (\beta_1,\id)^*\ev^*EG,\ \
\nabla' = \beta^*\nabla^u = (\beta_1,\id)^* (\ev^*\nabla^u).$$ Thus the
pair $(\ev^*EG, \ev^*\nabla^u)$ is a universal bundle with connection
for these bundles.  

We will need some topological converses. First, if $X\to B$ is a
$\calG$-bundle and if $\beta:B\times N\to BG$ classifies $X'$, then
$\beta_1$ classifies X; i.e.~$X\simeq \beta_1^*E\calG.$ The proof is
an easy check. Conversely, 
 if $X\to B$ is classified by a map $\beta_1:B\to
B\calG,$ then $(\beta_1,\id):B\times N\to B\calG\times N$ has
$(\beta_1,\id)^*E\calG' \simeq X'.$ Indeed,
\begin{eqnarray*} (\beta_1,\id)^*E\calG'|_{(b,x)} &=&
 E\calG'|_{(\beta_1(b),x)} =
\MMaps^{G}(Fr(F_x), EG|_{\beta_1(b)(x)})\\
&=&
{\rm Aut}(F_x,\gamma_n^\infty|_{\beta_1(b)(x)}) = \beta_1^*E\calG|_b
\ {\rm at}\  x.\end{eqnarray*}
This last expression equals $X|_b$ at $x$, which equals $X'|_{(b,x)}.$  

We now use the  fact that in this notation
$$X'\simeq (\beta_1,\id)^*E\calG'\simeq (\beta_1,\id)^*\ev^*EG.$$
to show that there exists a non-canonical isomorphism
$\alpha:E\calG' \to \ev^*EG.$  As a result, we get:

\begin{lemma} Let $X$ be a $\calG$ bundle over $B$, and let $X'$ 
be the associated $G$ bundle over $B\times N.$ Then the pair
$(E\calG',\underline{\alpha}^*\ev^*\nabla^u)$ is geometrically classifying for
$(X',\nabla')$, for any connection $\nabla'$ on $X'.$ \end{lemma}

Here $\underline{\alpha}:B\times N\to B\calG\times N$ is the map under 
 $\alpha.$

\begin{proof}
Let $a:B\calG \times N\to BG$ classify $E\calG',$ and let
$\beta:B\times N\to BG$ classify $X'$.  
Then as above $$X'\simeq (\beta_1,\id)^*\ev^*EG.$$ 
By the work above, we also have $X'\simeq (\beta_1,\id)^*E\calG'$ and
$E\calG' \simeq (a_1,\id)^*\ev^*EG$, so 
$$X'\simeq
(\beta_1,\id)^*(a_1,\id)^*\ev^*EG.$$
The last two indented equations imply that 
$$\ev\circ(\beta_1,\id) \sim
\ev\circ(a_1\beta_1,\id).$$

Now if $f_0,f_1:B\to B\calG$ have $\ev\circ (f_0, \id)\sim \ev\circ
(f_1,\id),$ then $f_0\sim f_1.$ For we have $F_t =
F(\cdot,\cdot,t):B\times N\times [0,1]\to BG$ with $F_0(b,x) =
\ev\circ(f_0,\id)(b,x) = f_0(b)(x)$, and $F_1(b,x) = f_1(b)(x).$ Define
$\bar f_t:B\to B\calG$ by $\bar f_t(b)(x) = F_t(b,x).$ Then $\bar f_0
= f_0, \bar f_1 = f_1.$

Thus $\beta_1\sim a_1\beta_1$ for all maps $\beta_1:B\to B\calG$.
Setting $\beta_1 = \id$ (which
corresponds to taking $\beta = \ev$), we get $\id\sim a_1.$  Thus 
$E\calG'\simeq (a_1,\id)^*\ev^* EG\simeq \ev^*EG$, and we set
$\alpha$
to be the composition of these isomorphisms. \end{proof}

By this lemma, we expect that $E\calG$ will also have a universal
connection for $\calG$ bundles with connection.  
 We need to be able to go back and forth between
connections on a $\calG$ bundle $X\to B$ and connections on the
associated $X'\to B\times N.$ One direction is clear: given a
connection $\nabla'$ on $X'$, we define a connection $\nabla = \nabla'_{\ext}$ 
on $X$
by
\begin{equation}\label{prime}\nabla_Z s (b)(x) = \nabla'_{(Z,0)} u_s(b,x),
\end{equation}
for $Z\in TB.$  

The only difficulty in going from $\nabla$ to $\nabla'$ is defining
differentiations in $N$ directions.  So given a connection $\nabla$ on $X$, we
use (\ref{prime}) to define $\nabla'_{(Z,0)}$, and for $Y\in TN$, we set
\begin{equation}\label{arb}\nabla'_{(0,Y)} u_s(b,x) = \nabla^b_Y u_s(b,x),
\end{equation}
where $\nabla^b$ is an arbitrary (but smoothly varying in $b$) connection on 
$X'|_{\{b\}\times N}.$  Note that $\nabla\mapsto\nabla'\mapsto\nabla$ is the
 identity, but $\nabla'\mapsto\nabla\mapsto\nabla'$ is not defined,
due to this arbitrariness.

As $\beta = (\beta_1,\id): B \times N\to B\calG$ runs through all maps
classifying $X'$ (equivalently, as $\beta_1:B\to B\calG$ runs through
all maps classifying $X$), and as $\alpha$ runs through all
isomorphisms of $X'$ and $\beta^* \ev^*EG$, we know that
$\alpha^{-1}\beta^* \ev^*\nabla^u$ runs through all connections on
$X'$.  Here $\alpha^{-1}\beta^* \ev^*\nabla^u u_s$ means
$\alpha^{-1}(\beta^* \ev^*\nabla^u (\alpha \circ u_s)).$  

We claim that $\widetilde \nabla \equiv
( \ev^*\nabla^u)_{\ext}$ is a universal connection on
$E\calG \to B\calG.$ Given a connection $\nabla$ on $X$, construct
$\nabla'$ as in (\ref{prime}), (\ref{arb}), and write
$\nabla' = \alpha^{-1}\beta^* \ev^*\nabla^u$ for some $\alpha,\beta.$
Then
$$\nabla_Z s(b)(x) = \nabla'_{(Z,0)}u_s(b,x) = 
\alpha^{-1}(\beta^* \ev^*\nabla^u_{(Z,0)} (\alpha \circ u_s))(b,x).$$
For $\overline{\beta_1}:
\beta_1^*E\calG \to E\calG$ 
the bundle map corresponding to $\beta_1,$ we have
\begin{eqnarray*} \alpha^{-1}(\beta_1^*\widetilde \nabla)_Z
(\alpha\circ s)(b)(x) &=& \alpha^{-1}\widetilde 
\nabla_{(\beta_1)_*Z} (\overline{\beta_1}\circ\alpha\circ s)(b)(x)\\
&=& \alpha^{-1}(\ev^*\nabla^u)
_{((\beta_1)_*Z,0)} u_{\overline{\beta_1}\circ\alpha\circ s}(b,x)\\
&=& \alpha^{-1}(\ev^*\nabla^u)_{\beta_*(Z,0)} u_{
\overline{\beta_1}\circ\alpha\circ s}(b,x) \\
&=& \alpha^{-1}(\beta^* \ev^*\nabla^u)_{(Z,0)} u_{\alpha\circ s}(b,x)\\
&=& \alpha^{-1}  
(\beta^* \ev^*\nabla^u)_{(Z,0)} (\alpha\circ u_s)(b,x).\end{eqnarray*}
In the last two lines, we have used 
$u_{\overline{\beta_1}\circ \alpha\circ s} =  \overline{\beta} \circ
u_{\alpha\circ s},\ \ 
u_{\alpha\circ s} = \alpha\circ u_s,$
for $\alpha\circ s \in\Gamma(B,\beta_1^*E\calG)$, which we leave as
 annoying exercises.  As a result, we have:

\begin{theorem} The natural extension of $\ev^*\nabla^u$ to a connection on 
$E\calG \to B\calG$ is a universal connection for $\calG$ bundles.
\end{theorem}
\bigskip

As a final remark, we show that the bundle  $E\calG'_0$ is isomorphic
to the universal bundle $\widetilde\PP$ in gauge theory \cite[Ch.~5]{DK}.
Here ${\mathcal G}_0$ is the 
group of gauge transformations
which are the identity at a fixed point of $N$, and
$\widetilde\PP = ({\mathcal A}\times P)/{\calG}$ as a
${\calG}$-bundle over $({\calAA}\times N)/{\calG} = B{\calG}\times
N$ with the action $g\cdot(\nabla,v) = (g\cdot \nabla, g_x^{-1}v)$ for
$v\in P_x.$  Recall that  $\calAA$ is the space of connections on the principal
bundle $P = \Fr(F).$

\begin{lemma} $ E\calG'_0 \simeq \widetilde\PP$  as bundles over $B{\mathcal
G}_0\times N.$  \end{lemma}

\begin{proof}   We have shown that $E{\mathcal G}'\simeq
\ev^*EG$, and the same argument works for ${\mathcal G}_0'$.

We drop the subscript $0.$
Since $E{\mathcal G} = \MMaps^{G}(P,EG)$, every element of
which covers some $f:N\to BG$, we can write
$$E{\mathcal G} = \{(\alpha,f): \alpha\in \Iso(P,f^*EG)\}.$$
There is a map $\Psi: E{\mathcal G}\to {\mathcal A}$, 
thought of as the space of all connections on $E$,
given by $\Psi(\alpha,f) = \alpha\cdot
(f^*\nabla^u) = \alpha^{-1}f^*\nabla^u\alpha.$  This map is surjective
by the finite dimensional universal connection theorem.  We define the
trivial equivalence relation $(\alpha,f)\sim (\beta,h) \Leftrightarrow
\Psi(\alpha,f) = \Psi(\beta,h).$  Thus the space of all connections
${\mathcal A}$ is in one-to-one correspondance with $E{\calG
}/\sim.$

From its definition,
$$\widetilde\PP = \frac{  \frac{E{\calG}}{(\alpha,f)\sim (\beta,h)}
\times P}{{\calG}},$$
as one checks that $(\alpha,f)\sim (\beta,h)$ iff $(\alpha\cdot g
,f)\sim (\beta\cdot g,h).$  Fiberwise, we have
$$\widetilde\PP|_{(f,x)} = \{ \left[ [\alpha_x,f],v\right]: \alpha_x \
{\rm is \ the\ restriction\ of}\ \alpha\in \Iso(P,f^*EG)\ {\rm to}\
x\},$$
where the equivalence class $\left[ [\alpha_x,f],v\right]$ is defined
by
$([\alpha_x,f],v) \sim ([\alpha_xg_x,f], g_x^{-1}v).$  

We define a map $\widetilde\PP\to \ev^*EG$ by $\left[
[\alpha_x,f],v\right]\mapsto \alpha_x(v)|_f.$  This map is well
defined and surjective.  For injectivity, say $\alpha_x(v) =
\alpha'_x(v').$  There exists $g$ with $\alpha'_x = \alpha_x g$, so
$\alpha'_x(v') = (\alpha_xg)(v') = \alpha_x(gv').$  Since $\alpha_x$
is an isomorphism, $v= gv'.$  Thus
$$\left[ [\alpha'_x,f],v'\right] = \left[[\alpha_xg,f],g^{-1}v\right]
=\left[ [\alpha_x,f],v\right].$$
\end{proof}

\bibliographystyle{plain}

\bibliography{pr2}

\begin{thebibliography}{10}

\bibitem{AB}
M.~F. Atiyah and R.~Bott.
\newblock The {Y}ang-{M}ills equations over {R}iemann surfaces.
\newblock {\em Philos. Trans. Roy. Soc. London Ser. A}, 308:523--615, 1983.

\bibitem{BGV}
N.~Berline, E.~Getzler, and M.~Vergne.
\newblock {\em Heat Kernels and Dirac Operators}.
\newblock Springer-Verlag, Berlin, 1992.

\bibitem{B}
J.-{M}. Bismut.
\newblock Localization formulae, super-connections, and the index theorem for
  families.
\newblock {\em Commun. Math. Phys}, 103:127--166, 1986.

\bibitem{BF}
J.-{M}. Bismut and D.~Freed.
\newblock The analysis of elliptic operators. {I}.
\newblock {\em Commun. Math. Phys}, 106:159--176, 1986.

\bibitem{CDMP}
A.~Cardona, C.~Ducourtioux, J.-{P}. Magnot, and S.~Paycha.
\newblock Weighted traces on algebras of pseudo-differential operators and
  geometry of loop groups.
\newblock Infinite Dim. Anal,. Quantum Prob. and Related Topics, to appear.

\bibitem{DK}
S.~Donaldson and P.~Kronheimer.
\newblock {\em The Geometry of Four-Manifolds}.
\newblock Oxford University Press, New York, 1990.

\bibitem{F}
D.~Freed.
\newblock The geometry of loop groups.
\newblock {\em J. Differential Geometry}, 28:223--276, 1988.

\bibitem{GS}
G.~Grubb and R.~Seeley.
\newblock Weakly parametric pseudodifferential opeartors and
  {A}tiyah-{P}atodi-{S}inger boundary problems.
\newblock {\em Inventiones Math.}, 121:481--529, 1995.

\bibitem{K}
C.~Kassel.
\newblock Le r\'esidu non commutatif [d'apr\'es {W}odzicki].
\newblock {\em Ast\'erisque}, 177-178:199--229, 1989.

\bibitem{Ku}
N.~Kuiper.
\newblock The homotopy type of the unitary group of {H}ilbert space.
\newblock {\em Topology}, 3:19--30, 1965.

\bibitem{jp}
J.-{P}. Magnot.
\newblock Sur la g\'eometrie d'espaces de lacets point\'es.
\newblock Th\`ese, Universit\'e Blaise Pascal (Clermont II), 2001.

\bibitem{NR}
M.~Narasimhan and S.~Ramanan.
\newblock Existence of universal connections.
\newblock {\em Amer. J. Math.}, 83:563--572, 1961.

\bibitem{PR}
S.~Paycha and S.~Rosenberg.
\newblock Curvature of determinant bundles and first {C}hern forms.
\newblock {\it J. Geometry and Physics}, to appear, 2002.

\bibitem{PR2}
S.~Paycha and S.~Rosenberg.
\newblock Traces and characteristic classes on loop spaces.
\newblock preprint, 2003.

\end{thebibliography}

\end{document}